\documentclass[english,11pt,reqno]{amsart}
\usepackage{amssymb,amsmath,amstext,amsthm,amsfonts}
\usepackage{graphicx}
\usepackage[ansinew]{inputenc}

\newcommand{\length}{\operatorname{length}}
\newcommand{\leb}{\operatorname{Leb}}
\newcommand{\loc}{\text{loc}}
\newcommand{\dist}{\operatorname{dist}}
\newcommand{\diam}{\operatorname{diam}}
\newcommand{\per}{\operatorname{Per}}

\begin{document}

\newcommand{\mcup}{\mbox{$\bigcup$}}
\newcommand{\mcap}{\mbox{$\bigcap$}}

\def \RR {{\mathbb R}}
\def \ZZ {{\mathbb Z}}
\def \NN {{\mathbb N}}
\def \PP {{\mathbb P}}
\def \TT {{\mathbb T}}
\def \II {{\mathbb I}}
\def \JJ {{\mathbb J}}

\def \vare {\varepsilon }

 \def \cf {\mathcal{F}}
 \def \cm {\mathcal{M}}
 \def \cn {\mathcal{N}}
 \def \cq {\mathcal{Q}}
 \def \cp {\mathcal{P}}
 \def \cb {\mathcal{B}}
 \def \cc {\mathcal{C}}
 \def \cs {\mathcal{S}}
 \def \bc {\mathcal{B}}
 \def \hc {\mathcal{C}}

\newcommand{\dem}{\begin{proof}}
\newcommand{\cqd}{\end{proof}}

\newcommand{\qand}{\quad\text{and}\quad}

\newtheorem{theorem}{Theorem}
\newtheorem{corollary}{Corollary}

\newtheorem*{Maintheorem}{Main Theorem}

\newtheorem{maintheorem}{Theorem}
\renewcommand{\themaintheorem}{\Alph{maintheorem}}
\newcommand{\cmt}{\begin{maintheorem}}
\newcommand{\fmt}{\end{maintheorem}}

\newtheorem{maincorollary}[maintheorem]{Corollary}
\renewcommand{\themaintheorem}{\Alph{maintheorem}}
\newcommand{\cmc}{\begin{maincorollary}}
\newcommand{\fmc}{\end{maincorollary}}

\newtheorem{T}{Theorem}[section]
\newcommand{\ct}{\begin{T}}
\newcommand{\ft}{\end{T}}

\newtheorem{Corollary}[T]{Corollary}
\newcommand{\cco}{\begin{Corollary}}
\newcommand{\fco}{\end{Corollary}}

\newtheorem{Proposition}[T]{Proposition}
\newcommand{\cpr}{\begin{Proposition}}
\newcommand{\fpr}{\end{Proposition}}

\newtheorem{Lemma}[T]{Lemma}
\newcommand{\cle}{\begin{Lemma}}
\newcommand{\fle}{\end{Lemma}}

\newtheorem{Sublemma}[T]{Sublemma}
\newcommand{\csle}{\begin{Lemma}}
\newcommand{\fsle}{\end{Lemma}}

\theoremstyle{definition}

\newtheorem{Example}[T]{Example}
\newcommand{\cex}{\begin{Example}}
\newcommand{\fex}{\end{Example}}

\newtheorem{Remark}[T]{Remark}
\newcommand{\cre}{\begin{Remark}}
\newcommand{\fre}{\end{Remark}}

\newtheorem{Definition}[T]{Definition}
\newcommand{\cd}{\begin{Definition}}
\newcommand{\fd}{\end{Definition}}

\title[Partially hyperbolic sets]
{Topological structure of (partially) hyperbolic sets with positive
volume}

\author{Jos\'e F. Alves}
\address{Jos\'e F. Alves\\ Departamento de Matem\'atica Pura, Faculdade de Ci\^encias do Porto\\
Rua do Campo Alegre 687, 4169-007 Porto, Portugal}
\email{jfalves@fc.up.pt} \urladdr{http://www.fc.up.pt/cmup/jfalves}

\author{Vilton  Pinheiro}
\address{Vilton  Pinheiro\\ Departamento de Matem\'atica, Universidade Federal da Bahia\\
Av. Ademar de Barros s/n, 40170-110 Salvador, Brazil.}
\email{viltonj@ufba.br}
%\urladdr{http://www.fc.up.pt/cmup/home/jfalves}

\date{\today}

\thanks{Work carried out at the Federal University of
Bahia, University of Porto and IMPA. JFA was partially supported by
CMUP, by a grant of FCT and by POCI/MAT/61237/2004. VP was partially
supported by PADCT/CNPq and by POCI/MAT/61237/2004}

\keywords{hyperbolic set, partially hyperbolic set, horseshoe}

\begin{abstract}
We consider both hyperbolic sets and partially hyperbolic sets
attracting a set of points with positive volume in a Riemannian
manifold. We obtain several results on the topological structure of
such sets for diffeomorphisms whose differentiability is bigger than
one. We show in particular  that there are no partially hyperbolic
horseshoes with positive volume for such diffeomorphisms. %This
%generalizes a classical result by Bowen for uniformly hyperbolic
%horseshoes.
We also give a description of the limit set of almost every point
belonging to a hyperbolic set or a partially hyperbolic set with
positive volume.

\end{abstract}

\maketitle

\setcounter{tocdepth}{1}

\tableofcontents

%\newpage

\section{Introduction}

Since the 60's that  hyperbolic sets have played an important role
in the development of
the Theory of Dynamical Systems. %, either continuous (vector fields in
%a manifold) or discrete (smooth maps from a manifold $M$ into
%itself).
These are  invariant (by a smooth map) compact sets over which the
tangent bundle splits into two invariant subbundles, one of them
contracting and the other one expanding under the action of the
derivative of the map.  In this work we are concerned with discrete
dynamical systems (smooth transformations of a manifold), but our
techniques proved also usefulness in the continuous setting (vector
fields in a manifold), specially for the study of
singular-hyperbolic sets done in \cite{AAPP}. In the last decades an
increasing emphasis is being put on the dynamics of partially
hyperbolic sets. These are compact invariant sets for which the
tangent bundle splits into two invariant subbundles having
contracting/expanding behavior in one
direction and the other one being dominated by it. %These will also
%be part of our study in this work.
Precise definitions of all these objects will be given in the next
section.

In this context a special role has been played by the horseshoes,
which have been introduced by Smale, and as shown in \cite{S},
always exist near a transverse homoclinic point associated to some
hyperbolic periodic point of saddle type, i.e. a point whose orbit
asymptotically approaches that saddle point, both in the past and in
the future. Horseshoes can be used to show that transverse
homoclinic points are always accumulated by periodic points, but the
dynamical richness of these objects goes far beyond the initial
application by Smale, and since then many other results have been
proved using horseshoes. These are Cantor sets which are, in
dynamical terms, topologically conjugated to full shifts.

A special interest lies in the horseshoes that appear when one
unfolds a homoclinic tangency. Knowing how {\em fat} these
horseshoes are can have several implications in the dynamical
behavior after the homoclinic bifurcation. In this setting we
mention the {\em thickness}, which has been used by Newhouse
\cite{N} to prove the existence of infinitely many sinks, and the
{\em Hausdorff dimension}, which has been used by Palis, Takens and
Yoccoz to study the prevalence of hyperbolicity after the unfolding
of a homoclinic tangency; see \cite{PT,PY,PY2}.

One interesting issue we will be addressed to is the {\em volume} of
horseshoes. As shown by Bowen in \cite{Bo}, there are $C^1$
diffeomorphisms with hyperbolic horseshoes of positive volume. On
the other hand, Bowen has proved in \cite[Theorem 4.11]{B} that a
basic set (locally maximal hyperbolic set with a dense orbit) of a
$C^2$ diffeomorphism which attracts a set with positive volume,
necessarily attracts a neighborhood of itself. In particular, the
unstable manifolds through points of this set must be contained in
it, and consequently $C^2$ diffeomorphisms have no horseshoes with
positive volume.

For diffeomorphisms whose differentiability is higher than one, we
prove the nonexistence of horseshoes with positive volume in a more
general context of sets with some partially hyperbolic structure.
Using our framework in the context of hyperbolic sets, we are able
to show that Bowen's result  still holds without the local
maximality assumption, i.e. a transitive hyperbolic set which
attracts a set with positive volume necessarily attracts a
neighborhood of itself.
 Furthermore, we are able to prove
that there are no proper transitive hyperbolic sets with positive
volume for diffeomorphisms whose differentiability is higher than
one. Similar results for sets with nonempty interior had already
been obtained in \cite[Theorem 1]{ABD} and in \cite[Theorem 1.1]{F}.
%The proof of this result is based on the Thermodynamical formalism.
%We make a different proof and get rid of the locally maximal
%assumption. In particular, such a transitive hyperbolic set
On the other hand, as described in \cite[Remark 2.1]{ABD} or in
\cite[Example 2]{F}, there exist (non-transitive) hyperbolic
sets with positive volume which do not attract neighborhoods of themselves. %which %do not coincide with the whole manifold
%thus showing that  the transitivity.

 Let us mention two more important results in
this direction. It follows from \cite[Theorem~2]{Y}  that proper
uniformly partially hyperbolic sets supporting a unique equilibrium
state and attracting open neighborhoods of themselves necessarily
have zero volume. In the conservative setting, \cite[Theorem 15]{BV}
gives that a hyperbolic set for a volume preserving $C^2$
diffeomorphism either has zero volume or coincides with the whole
manifold.

In this work we also give a good description of the limit set of
almost every point in a hyperbolic set with positive volume: there
is a finite number of basic sets for which the $\omega$-limit set of
Lebesgue almost every point is contained in one of these basic sets.
We are also able to prove in a partially hyperbolic setting that
these $\omega$-limit sets are contained in the closure of finitely
many hyperbolic periodic points.

\subsection*{Acknowledgement} We are grateful to M. Viana for
valuable discussions and references on these topics.

\section{Statement of results}

Let $f:M\to M$ be a diffeomorphism of a compact connected Riemannian
manifold $M$. We say that $f$ is $C^{1+}$ if $f$ is $C^1$ and $Df$
is H\"older continuous. We use $\leb$ to denote a normalized volume
form on the Borel sets of $M$ that we call {\em Lebesgue measure.}
Given a submanifold $\gamma\subset M$ we use $\leb_\gamma$ to denote
the measure on $\gamma$ induced by the restriction of the Riemannian
structure to $\gamma$. A set  $\Lambda\subset M$ is said to be {\em
invariant} if $f(\Lambda)=\Lambda$, and {\em positively invariant}
if $f(\Lambda)\subset\Lambda$.

%In the conservative case, Xia proved that if a
%$C^{1+}$ diffeomorphism has a hyperbolic set with positive
%Lebesgue measure, then $f$ must necessarily be an Anosov
%diffeomorphism, i.e. the whole manifold $M$ is a hyperbolic set.
%

\subsection{Partially hyperbolic sets}
Let $K$ be a positively invariant compact  set, and define
 $$\Lambda=\bigcap_{n\ge0}f^n(K).$$
Suppose  that there is a continuous splitting $T_K M=E^{cs}\oplus
E^{cu}$ of the tangent bundle restricted to $K$, and assume that
this splitting is $Df$-invariant over~$\Lambda$. We say that this is
a  {\em dominated splitting\/} (over $\Lambda$) if there is a
constant $0<\lambda<1$ such that for some choice of a Riemannian
metric on $M$
$$
 \|Df \mid E^{cs}_x\|
\cdot \|Df^{-1} \mid E^{cu}_{f(x)}\| \le\lambda,\quad \text{for
every $x\in\Lambda$}.
$$
 We call $E^{cs}$ the {\em centre-stable
bundle} and $E^{cu}$ the {\em centre-unstable bundle}. $\Lambda$~is
said to be {\em partially hyperbolic}, if additionally $E^{cs}$ is
{\em uniformly contracting} or $E^{cu}$ is {\em uniformly
expanding}, meaning that there exists
%$$
%\|Df\mid E^{cs}\| \le\lambda\quad\text{or}\quad \|Df^{-1}\mid
%E^{cu}\| \le\lambda.
%$$
$0<\lambda<1$ such that
 $$ \|Df\mid E^{cs}_x\|
\le\lambda, \quad \text{for every $x\in\Lambda$},
 $$
 or
 $$\|Df^{-1}\mid E^{cu}_{f(x)}\| \le\lambda,\quad \text{for
every $x\in\Lambda$}.
 $$
%for every $x\in \Lambda$. %We shall write $E^{s}$ instead of $E^{cs}$
%in the first case, and $E^{s}$ instead of $E^{cu}$ in the other
%case.
 We say that $f$ is {\em non-uniformly expanding along the
centre-unstable direction} in $K$ if there is $c>0$ such that for
Lebesgue almost every $x\in K$
\begin{equation}
\label{NUE}\tag{NUE} \liminf_{n\to+\infty} \frac{1}{n}
    \sum_{j=1}^{n} \log \|Df^{-1} \mid E^{cu}_{f^j(x)}\|<-c.
\end{equation}
Condition \ref{NUE} means that the derivative  has {\em expanding
behavior in the centre-unstable direction in average} over the orbit of $x$ for an infinite number of times. %This implies that $x$
%has $\dim(E^{cu})$ positive
%Lypaunov exponents in the $E^{cu}_x$ direction. %In particular,
%local unstable manifolds exist for points in $\Lambda$,
%by~\cite{Pe76}.
If condition~\ref{NUE} holds for every point in a compact invariant
set   $\Lambda$, then $E^{cu}$ is {\em uniformly expanding in the
centre-unstable direction in $\Lambda$}. This is not necessarily the
case if \ref{NUE} occurs only Lebesgue almost everywhere. A class of
diffeomorphisms with a dominated splitting $TM=E^{cs}\oplus E^{cu}$
for which \ref{NUE} holds Lebesgue almost everywhere in $M$ and
$E^{cu}$ is not uniformly expanding can be found in \cite[Appendix
A]{ABV}.

%Horseshoes were introduced by Smale and appear naturally when one
%unfolds a homoclinic tangency associated to some hyperbolic periodic
%point of saddle type. It follows from \cite[Theorem 4.11]{B}  that a
%$C^{1+}$ diffeomorphism cannot have a fat hyperbolic horseshoe, i.e.
%a hyperbolic horseshoe $\Lambda$ with $\leb(\Lambda)>0$; actually
%the result in~\cite{B} is proved for basic sets. Let us remark that
%fat hyperbolic horseshoes exist for $C^1$ diffeomorphisms, as shown
%in \cite{Bo}.

We say that an embedded disk $\gamma\subset M$ is an {\em unstable
manifold}, or an {\em unstable disk}, if
$\dist(f^{-n}(x),f^{-n}(y))\to0$ exponentially fast as $n\to\infty$,
for every $x,y\in\gamma$. Similarly, $\gamma$ is called a {\em
stable manifold}, or a {\em stable disk}, if
$\dist(f^{n}(x),f^{n}(y))\to0$ exponentially fast as $n\to\infty$,
for every $x,y\in\gamma$. It is well-known that every point in a
hyperbolic set possesses a local stable manifold $W_{loc}^s(x)$ and
a local unstable manifold $W_{loc}^u(x)$ which are disks tangent to
$E^{s}_x$ and $E^{u}_x$ at $x$ respectively. A compact invariant set
$\Lambda$ is said to be {\em horseshoe-like} if there are local
stable and local unstable manifolds through all its points which
intersect $\Lambda$ in a Cantor set.

 \cmt\label{t:disco1}
 Let \( f: M\to M \) be a  \( C^{1+} \)
diffeomorphism and let $K\subset M$ be a forward invariant compact
 set with a continuous splitting $T_K M=E^{cs}\oplus E^{cu}$
dominated over $\Lambda=\bigcap_{n\ge0}f^n(K)$. If \ref{NUE} holds
for a positive Lebesgue set of points $x\in K$, then $\Lambda$
contains some local unstable disk.
 \fmt

 The next result is a direct consequence of Theorem~\ref{t:disco1},
whenever $E^{cu}$ is uniformly expanding. If, on the other hand,
$E^{cs}$ is uniformly contracting, then we just have to apply
Theorem~\ref{t:disco1} to $f^{-1}$.

 \cmc\label{c:disco1}
  Let \( f: M\to M \) be a  \( C^{1+} \)
diffeomorphism and let $K\subset M$ be a compact invariant set with
$\leb(K)>0$  having a continuous splitting $T_K M=E^{cs}\oplus
E^{cu}$ for which $\Lambda=\bigcap_{n\ge0}f^n(K)$ is partially
hyperbolic.
\begin{enumerate}
  \item If $E^{cs}$ is uniformly contracting, then $\Lambda$ contains a local stable disk.
  \item If $E^{cu}$ is uniformly expanding, then $\Lambda$ contains a local unstable disk.
\end{enumerate}

 \fmc

In particular, \( C^{1+} \) diffeomorphisms have no partially
hyperbolic horseshoe-like sets with positive volume. The same
conclusion holds for partially hyperbolic sets intersecting a local
stable disk or a local unstable disk in a positive Lebesgue measure
subset, as Corollary~\ref{c:disco2} below shows.

 \cmt\label{t:disco2}
 Let \( f: M\to M \) be a  \( C^{1+} \)
diffeomorphism and let $K\subset M$ be a forward invariant compact
set with a continuous splitting $T_K M=E^{cs}\oplus E^{cu}$
dominated over $\Lambda=\bigcap_{n\ge0}f^n(K)$.  Assume that there
is a local unstable disk $\gamma$
 such that \ref{NUE} holds for
every $x$ in a positive $\leb_\gamma$ subset of $\gamma\cap K$. Then
$\Lambda$ contains some local unstable disk.
 \fmt

The next result is an immediate consequence of
Theorem~\ref{t:disco2}, in the case that $E^{cu}$ is uniformly
expanding, and a consequence of the same theorem applied to $f^{-1}$
when $E^{cs}$ is uniformly contracting. Actually, we shall prove a
stronger version of this result in Theorem~\ref{t:disco}.

 \cmc\label{c:disco2}
 Let \( f: M\to M \) be a  \( C^{1+} \)
diffeomorphism and let $K\subset M$ be a forward invariant compact
set having a continuous splitting $T_K M=E^{cs}\oplus E^{cu}$
dominated over $\Lambda=\bigcap_{n\ge0}f^n(K)$.
\begin{enumerate}
\item If $E^{cs}$ is uniformly contracting and there is a local stable disk
$\gamma$ such that $\leb_\gamma(\gamma\cap K)>0$, then $\Lambda$
contains a local stable disk.
\item If $E^{cu}$ is
uniformly expanding and there is a local unstable disk $\gamma$ such
that $\leb_\gamma(\gamma\cap K)>0$, then $\Lambda$ contains a local
unstable disk.
\end{enumerate}
 \fmc

Using the previous results we are able to give a description of the
$\omega$-limit of Lebesgue almost every point in  a partially
hyperbolic whose center-unstable direction displays non-uniform
expansion in a subset with positive volume. Recall that the {\em
$\omega$-limit} of  $x\in M$ is the set of accumulation points of
its orbit.

\cmt\label{t:limitph}
 Let \( f: M\to M \) be a  \( C^{1+} \)
diffeomorphism and let $K\subset M$  be a forward invariant compact
set with $\leb(K)>0$  having a continuous splitting $T_K
M=E^{cs}\oplus E^{cu}$ for which $\Lambda=\bigcap_{n\ge0}f^n(K)$ is
partially hyperbolic. Assume that $E^{cs}$ is uniformly contracting
and \ref{NUE} holds for Lebesgue almost every $x\in K$. Then there
are hyperbolic periodic points $p_1,\dots,p_k\in\Lambda$ such that:
\begin{enumerate}
    \item $ \overline{W^u(p_i)}\subset \Lambda$ for each $1\le i\le
    k$;
    \item for $\leb$ almost every $x\in K$ there is $1\le i\le
    k$ with
    $\omega(x)\subset \overline{W^u(p_i)}$.
    %\item $ \overline{W^u(p_i)}$ attracts an open set $A_i$ containing $ W^u(p_i)$ for each $1\le i\le
%    k$.
\end{enumerate}
Moreover, if $E^{cu}$ has dimension one, then for each $1\le i\le
    k$
 \begin{enumerate}
    \item[(3)] $\overline{W^u(p_i)}$ attracts an open neighborhood of itself. %set containing $ \overline{W^u(p_i)}$

\end{enumerate}
 \fmt

This last conclusion also holds whenever $E^{cs}$ is uniformly
contracting. Actually, more can be said in the case of uniformly
hyperbolic sets with positive volume as we show in the next
subsection.

\subsection{Hyperbolic sets}
 We say that a compact invariant set $\Lambda$ is {\em hyperbolic}
if there is a $Df$-invariant splitting $T_\Lambda M=E^{s}\oplus
E^{u}$ of the tangent bundle restricted to~$\Lambda$ and a constant
$\lambda<1$ such that (for some choice of a Riemannian metric on
$M$) for every $x\in \Lambda$
      $$\|Df \mid E^{s}_x\|<\lambda\qand \|Df^{-1} \mid E^{u}_x\|<\lambda.$$

We are able to prove that transitive hyperbolic sets with positive
volume necessarily coincide with the whole manifold, i.e. the
diffeomorphism is Anosov.

\cmt\label{t:Anosov} Let \( f: M\to M \) be a  \( C^{1+} \)
diffeomorphism and let $\Lambda\subset M$ be a transitive hyperbolic
set. \begin{enumerate}
    \item If $\Lambda$ has  positive volume, then $\Lambda=M$.
    \item If $\Lambda$ attracts a set with positive volume, then $\Lambda$ attracts a neighborhood of
    itself.
\end{enumerate}
 \fmt

The main reason why we cannot generalize the results in this
subsection to the context of partially hyperbolic sets is that the
length of local stable/unstable manifolds may shrink to zero when
iterated back/forth, respectively. The next result gives a
description of the $\omega$-limit of Lebesgue almost every point in
a hyperbolic set with positive volume. Taking $f^{-1}$ a similar
decomposition holds for $\alpha$-limits.

\cmt[Spectral Decomposition]\label{t:limit}
 Let \( f: M\to M \) be a  \( C^{1+} \)
diffeomorphism and let $\Lambda\subset M$ be a hyperbolic set with
positive volume.  There are hyperbolic sets
 $\Omega_1,\dots,\Omega_q\subset\Lambda$  such that:
\begin{enumerate}
    \item for $\leb$ almost every $x\in\Lambda$ there is $1\le i\le q$
    such that $\omega(x)\subset \Omega_i$;
    \item $\Omega_j$
    attracts a neighborhood of itself in $M$, for each $1\le j\le
    q$;
    %\item each $f\vert\mho_i$  and $f\vert\Omega_i$ is transitive
%    and contains dense periodic points.
    \item  $f\vert\Omega_k$ is transitive;
    \item $\per(f)$ is dense in $\Omega_j$,
    for each $1\le j\le q$.
\end{enumerate}
    Moreover, for each $1\le k\le q$ there is a decomposition of
    $\Omega_k$ into disjoint hyperbolic sets
    $\Omega_k=\Omega_{k,1}\cup\cdots\cup\Omega_{k,n_k}$ such that:
\begin{enumerate}
      \item[(5)] $f(\Omega_{k,i})=\Omega_{k,i+1}$, for $1\le i< n_k$,
      and $f(\Omega_{k,n_k})=\Omega_{k,1}$;
      \item[(6)] $f^{n_k}\colon \Omega_{k,i}\to\Omega_{k,i}$ is
      topologically mixing for every $1\le i\le n_k$.
    \end{enumerate}
    %An analogous decomposition holds for $\mho_k$, for each $1\le k\le p$.
 \fmt

\subsection{Overview} This paper is organized in the following way. In
Section~\ref{s.holder} we present some results from \cite{ABV} on
the H\"older control of the tangent direction of certain
submanifolds, and %then, using non-uniform expansion and hyperbolic
%times,
in Section~\ref{s.hyptimes} we derive some bounded distortion
results. Theorem~\ref{t:disco1} and Theorem~\ref{t:disco2} are
actually corollaries of a slightly more general result that we
present at the beginning of Section~\ref{s.local}. Let us mention
that the results in Section \ref{s.local} (specially
Lemma~\ref{l.discao}) are not consequence of the results in
\cite{ABV}, since we are using a weaker form of non-uniform
expansion in \ref{NUE}. Theorem~\ref{t:limitph} is proved in
Section~\ref{s.hypeper}. Finally, in Section~\ref{s.hyperbolic} we
%give an example of a $C^{1^+}$ diffeomorphism with a positive volume
%hyperbolic set not coinciding with the whole manifold, and
prove
Theorem~\ref{t:Anosov} and Theorem~\ref{t:limit}. %We finally  prove
%Theorem~\ref{t:Anosovconserva} in Section~\ref{s.PH}.

\section{H\"older control of tangent direction}\label{s.holder}
In this we present some results in~\cite[Section 2]{ABV} concerning
the H\"older control of the tangent direction of submanifolds.
Although those results are stated for $C^2$ diffeomorphims, they are
valid for diffeomorphisms of class $C^{1+}$, as observed in
\cite[Remark 2.3]{ABV}.

Let $K$ be a positively invariant compact  set for which there is a
continuous splitting $T_K M=E^{cs}\oplus E^{cu}$ of the tangent
bundle restricted to $K$ which is $Df$-invariant over
 $$\Lambda=\bigcap_{n\ge0}f^n(K).$$
 We fix continuous extensions of
the two bundles $E^{cs}$ and $E^{cu}$ to some compact neighborhood
$U$ of $\Lambda$, that we still denote by ${E}^{cs}$ and ${E}^{cu}$.
Replacing~$K$ by a forward iterate of it, if necessary, we may
assume that $K\subset U$.

Given $0<a<1$, we define the {\em centre-unstable cone field
$\left(C_a^{cu}(x)\right)_{x\in U}$ of width $a$\/} by
\begin{equation}
\label{e.cucone} C_a^{cu}(x)=\big\{v_1+v_2 \in {E}_x^{cs}\oplus
{E}_x^{cu} \text{\ such\ that\ } \|v_1\| \le a \|v_2\|\big\}.
\end{equation}
We define the {\em centre-stable cone field
$\left(C_a^{cs}(x)\right)_{x\in U}$ of width $a$\/} in a similar
way, just reversing the roles of the subbundles in (\ref{e.cucone}).

We fix $a>0$ and $U$ small enough so that, up to slightly increasing
$\lambda<1$, the domination condition remains valid for any pair of
vectors in the two cone fields, i.e.
$$
\|Df(x)v^{cs}\|\cdot\|Df^{-1}(f(x))v^{cu}\|
\le\lambda\|v^{cs}\|\,\|v^{cu}\|,
$$
for every $v^{cs}\in C_a^{cs}(x)$, $v^{cu}\in C_a^{cu}(f(x))$, and
any  $x\in U\cap f^{-1}(U)$. Note that the centre-unstable cone
field is positively invariant: $$Df(x) C_a^{cu}(x)\subset
C_a^{cu}(f(x)),\quad\text{whenever $x,f(x)\in U$.}$$ Indeed, the
domination property together with the invariance of $E^{cu}$ over
$\Lambda$ imply that
\begin{equation}\label{eqcone}
Df(x) C_a^{cu}(x) \subset C_{\lambda a}^{cu}(f(x))
                  %\subset C_a^{cu}(f(x)),
\end{equation}
 for every $x\in \Lambda$. This extends to any $x\in
U\cap f^{-1}(U)$ just by continuity, slightly increasing
$\lambda<1$, if necessary.

\cre The invariance of the splitting $T_K M=E^{cs}\oplus E^{cu}$  is
used in~\cite{ABV} to derive conclusions for the points in the small
neighborhood $U$ of $K$. Although here we are taking the invariance
of the splitting just restricted to $\Lambda$, since we are assuming
$K\subset U$, where $U$ is a small neighborhood of $\Lambda$, the
results of \cite[Section 2.1]{ABV} are still valid in our situation.
See also \cite[Remark 2.1]{ABV}. \fre

%Wherever we presume $E^{cs}$ to be uniformly contracting (as
%already mentioned, this will happen not happen before
%Subsection~\ref{ss.SRB}), we denote it by $E^{s}$ instead, and
%represent by $\tilde{E}^{s}$ its extension to $U$. Moreover, in
%that case the cone field $C^{cs}_a$ is denoted $C_a^{s}$, and
%called {\em strong-stable\/}.

%\subsection{Unstable direction}

We say that an embedded $C^1$ submanifold $N\subset U$ is {\em
tangent to the centre-unstable cone field\/} if the tangent subspace
to $N$ at each point $x\in N$ is contained in the corresponding cone
$C_a^{cu}(x)$. Then $f(N)$ is also tangent to the centre-unstable
cone field, if it is contained in $U$, by the domination property.
%The tangent bundle of $N$ is said to be H\"older continuous if $x
%\mapsto T_x N$ defines a H\"older continuous section from $N$ to the
%corresponding Grassman bundle of $M$. %In this subsection we show
%%that the tangent bundle of the iterates of a $C^2$ submanifold are
%%H\"older continuous (as long as they do not leave $U$ ), with
%%uniform H\"older constants.

 We choose $\delta_0>0$ small
enough so that the inverse of the exponential map $\exp_x$ is
defined on the $\delta_0$ neighbourhood of every point $x$ in $U$.
From now on we identify this neighbourhood of $x$ with the
corresponding neighbourhood $U_x$ of the origin in $T_x N$, through
the local chart defined by $\exp_x^{-1}$. %Accordingly, we identify
%$x$ with $0\in T_xN$.
Reducing $\delta_0$, if necessary, we may suppose that ${E}^{cs}_x$
is contained in the centre-stable cone $C^{cs}_a(y)$ of every $y\in
U_x$. In particular, the intersection of $C^{cu}_a(y)$ with
${E}^{cs}_x$ reduces to the zero vector. Then, the tangent space to
$N$ at $y$ is parallel to the graph of a unique linear map
$A_x(y):T_x N \to {E}_x^{cs}$. Given constants $C>0$ and $0<\zeta\le
1$, we say that {\em the tangent bundle to $N$ is
$(C,\zeta)$-H\"older\/} if for every $y\in N \cap U_x $ and $x\in
V_0$
\begin{equation}
\label{e.holder} \|A_x(y)\|\le C d_x(y)^\zeta ,
\end{equation}
where $d_x(y)$ denotes the distance from $x$ to $y$ along $N\cap
U_x$, defined as the length of the shortest curve connecting $x$ to
$y$ inside $N\cap U_x$.

Recall that we have chosen the neighbourhood $U$ and the cone width
$a$ sufficiently small so that the domination property remains valid
for vectors in the cones $C_a^{cs}(z)$, $C_a^{cu}(z)$, and for any
point $z$ in $U$. Then, there exist $\lambda_1 \in (\lambda,1)$ and
$\zeta\in(0,1]$ such that
\begin{equation}
\label{e.dominacao} \|Df(z) v^{cs}\| \cdot \|Df^{-1}(f(z))
v^{cu}\|^{1+\zeta} \le \lambda_1 < 1
\end{equation}
for every norm $1$ vectors $v^{cs}\in C_a^{cs}(z)$ and $v^{cu}\in
C_a^{cu}(z)$, at any $z\in U$. Then, up to reducing $\delta_0>0$ and
slightly increasing $\lambda_1<1$, condition (\ref{e.dominacao})
remains true if we replace $z$ by any $y\in U_x$, with $x\in U$
(taking $\|\cdot\|$ to mean the Riemannian metric in the
corresponding local chart).

We fix $\zeta$ and $\lambda_1$ as above. Given a $C^1$ submanifold
$N\subset U$, we define
\begin{equation}
\label{e.kappa} \kappa(N)=\inf\{C>0:\text{the tangent bundle of $N$
is $(C,\zeta)$-H\"older}\}.
\end{equation}
The next result appears in \cite[Corollary 2.4]{ABV}.

%
%\cpr \label{p.curvature} There exist $\lambda_0<1$ and $C_0>0$ so
%that if $N\subset U \cap f^{-1}(U)$ is any $C^1$ submanifold tangent
%to the centre-unstable cone field, then
%$$
%\kappa(f(N)) \le \lambda_0 \, \kappa(N) + C_0\,.
%$$
%\fpr
%

\cpr \label{c.curvature} There exists $C_1>0$ such that, given any
$C^1$ submanifold $N\subset U$ tangent to the centre-unstable cone
field,
\begin{enumerate}
\item there exists $n_0\ge 1$ such that $\kappa(f^n(N)) \le
C_1$ for every $n\ge n_0$ such that $f^k(N) \subset U$ for all $0\le
k \le n$;
 \item if $\kappa(N) \le C_1$, then the same is true
for every iterate $f^n(N)$  such that $f^k(N)\subset U$ for all
$0\le k \le n$;
 \item in particular, if $N$ and $n$ are
as in (2), then the functions
$$
J_k: f^k(N)\ni x \longmapsto \log |\det \big(Df \mid T_x
f^k(N)\big)|, \quad\text{$0\le k \le n$},
$$
are $(L,\zeta)$-H\"older continuous with $L>0$ depending only on
$C_1$ and~$f$.
\end{enumerate}
\fpr

\section{Hyperbolic times and bounded distortion}\label{s.hyptimes}

Let $K\subset M$ be a forward invariant compact set and let
$\Lambda\subset K\subset U$ be as in Section~\ref{s.holder}. The
following notion will allow us to derive {\em uniform behaviour\/}
(expansion, distortion) from the non-uniform expansion.

\cd \label{d.hyperbolic1} Given $\sigma<1$, we say that $n$ is a
{\em $\sigma$-hyperbolic time\/} for  $x\in K$ if
$$
\prod_{j=n-k+1}^{n}\|Df^{-1} \mid E^{cu}_{f^j(x)}\| \le \sigma^k,
\quad\text{for all $1\le k \le n$.}
$$
\fd %In particular, if $n$ is a $\sigma$-hyperbolic time for $x$,
%then $Df^{-k} \mid E^{cu}_{f^{n}(x)}$ is a contraction for every
%$1\le k \le n$:
%$$
%\|Df^{-k} \mid E^{cu}_{f^{n}(x)}\| \le \prod_{j=n-k+1}^{n}\|Df^{-1}
%\mid E^{cu}_{f^{j}(x)}\| \le \sigma^{k}.
%$$
If $a>0$ is taken sufficiently small in the definition of our cone
fields, and we choose $\delta_1>0$ also small so that the
$\delta_1$-neighborhood of $K$ should be contained in $U$, then by
continuity
\begin{equation}
\label{e.delta1} \|Df^{-1}(f(y)) v \| \le \frac{1}{\sqrt\sigma}
\|Df^{-1}|E^{cu}_{f(x)}\|\,\|v\|,
\end{equation}
whenever $x\in K$, $\dist(x,y)\le \delta_1$ and $v\in
C^{cu}_a(f(y))$.

%Let $D$ be a centre-unstable disk for which \eqref{NUE} holds
%$\leb_D$ almost everywhere.

Given any disk $\Delta\subset M$, we use $\dist_\Delta(x,y)$ to
denote the distance between $x,y\in \Delta$ measured along $\Delta$.
The distance from a point $x\in \Delta$ to the boundary of $\Delta$
is $\dist_\Delta(x,\partial \Delta)= \inf_{y\in\partial
\Delta}\dist_\Delta(x,y)$.

\smallskip

\cle \label{l.contraction} Take any $C^1$ disk $\Delta\subset U$ of
radius $\delta$, with $0<\delta<\delta_1$, tangent to the
centre-unstable cone field.
%Let $\Delta$ be a $C^1$ centre-unstable disk contained in $U$.
There is $n_0\ge1$ such that for $x\in \Delta\cap K$ with
$\dist_\Delta(x,\partial \Delta)\ge \delta/2$ and $n \ge n_0$  is a
$\sigma$-hyperbolic time for $x$, then there is a neighborhood $V_n$
of $x$ in $\Delta$ such that:
\begin{enumerate}
    \item $f^{n}$ maps $V_n$ diffeomorphically onto a
disk  of radius $\delta_1$ around  $f^{n}(x)$ tangent to the
centre-unstable cone field;
    \item for every $1\le k
\le n$ and $y, z\in V_n$, $$
\dist_{f^{n-k}(V_n)}(f^{n-k}(y),f^{n-k}(z)) \le
\sigma^{k/2}\dist_{f^n(V_n)}(f^{n}(y),f^{n}(z));$$

\item for every $1\le k
\le n$ and $y\in V_n$,
$$
\prod_{j=n-k+1}^{n}\|Df^{-1} \mid E^{cu}_{f^{j}(y)}\|
   \le \sigma^{k/2}.
$$
% \item if moreover $\kappa(\Delta) \le C_1$, then there
%is $C_2>1$ such that for $y,z\in V_n$
%$$
%\frac{1}{C_2} \le \frac{|\det Df^{n} \mid T_y \Delta|}
%                     {|\det Df^{n} \mid T_z \Delta|}
%            \le C_2.
%$$
\end{enumerate}
\fle

\begin{proof}
First we show that $f^n(\Delta)$ contains some disk of radius
$\delta_1$ around $f^n(x)$, as long as
 \begin{equation}\label{e.ene}
 n>2\frac{\log(\delta/(2\delta_1))}{\log(\sigma)}.
 \end{equation}
Define $\Delta_1$ as the connected component of $f(\Delta)\cap U$
containing $f(x)$. For $k\ge 1$, we inductively define
$\Delta_{k+1}\subset f^{k+1}(\Delta)$ as the connected component of
$f(\Delta_{k})\cap U$ containing $f^{k+1}(x)$. We shall prove that
$\Delta_n$ contains some disk of radius $\delta_1$ around $f^n(x)$,
for $n$ as in \eqref{e.ene}. Observe that since $\Delta_j\subset U$,
the invariance \eqref{eqcone} gives that for every $j\ge 1$
\begin{equation}\label{conetido}
T_w\Delta_j\subset C_{\lambda^j a}^{cu}(w),\quad\text{for every
$w\in\Delta_j$.}
\end{equation}

%Let $\vare=\min\{delta_1, \dist($ be the maximal radius for which
%the ball of radius $\vare$ inside $\Delta_n$ is contained
%$\dist_{\Delta_n}(f^n(x),f^n(y))<\delta_1$.

Let $\eta_0$ be a curve of minimal length in $\Delta_n$ connecting
$f^n(x)$ to $f^n(y)\in \Delta_n$ for which
$\dist_{\Delta_n}(f^n(x),f^n(y))<\delta_1$. For $0\le k \le n$,
writing $\eta_k=f^{-k}(\eta_0)$ we have $\eta_k\subset
\Delta_{n-k}$. We prove by induction that $\length(\eta_k) <
\sigma^{k/2}\delta_1, $ for $0 \le k \le n$. Let $1 \le k \le n$ and
assume that
$$
\length(\eta_j) < \sigma^{j/2}\delta_1, \quad\text{for\ }0 \le j \le
k-1.
$$
Denote by $\dot\eta_0(w)$ the tangent vector to the curve $\eta_0$
at the point $w$. Using the fact that $\eta_k\subset \Delta_{n-k}$
and
 \eqref{conetido} we have $$Df^{-j}(w)\dot\eta_0(w)\in
C_{\lambda^{n-j} a}^{cu}(f^{-j}(w))\subset C_{a}^{cu}(f^{-j}(w)).
$$ Then, by the choice of $\delta_1$ in (\ref{e.delta1}) and the
definition of $\sigma$-hyperbolic time,
$$
\|Df^{-k}(w) \dot\eta_0(w)\| \le {\sigma^{-k/2} \,
\|\dot\eta_0(w)\|}
         {\prod_{j=n-k+1}^{n}\|Df^{-1} |E^{cu}_{f^j(x)}\|}
\le \sigma^{k/2}\|\dot\eta_0(w)\|.
$$
Hence,
$$
\length(\eta_k) \le \sigma^{k/2}\length(\eta_0) < \sigma^{k/2}
\delta_1\,.
$$
This completes our induction.

In particular we have $ \length(\eta_n)  < \sigma^{n/2} \delta_1\,.
$ Moreover, the $k$ preimage of the ball of radius $\delta_1$ in
$\Delta_n$  centered at $f^n(x)$ is contained in $U$ for each $1\le
k\le n$.
 If $\eta_n$ is
a curve in $\Delta$ connecting $x$ to $y\in
\partial \Delta$, then we must have
 \begin{equation*}\label{e.ene2}
 n<2\frac{\log(\delta/(2\delta_1))}{\log(\sigma)}.
 \end{equation*}
 Hence  $f^n(\Delta)$ contains some disk of radius
$\delta_1$ around $f^n(x)$ for $n$ as in~\eqref{e.ene}.

Let now $D_1$ be the disk of radius $\delta_1$ around $f^n(x)$ in
$f^n(\Delta)$ and let $V_n=f^{-n}(D_1)$,  for $n$ as in
\eqref{e.ene}. Take any $y,z\in V_n$ and let $\eta_0$ be a curve of
minimal length in $D_1$ connecting $f^n(y)$ to $f^n(z)$. Defining
$\eta_k=f^{n-k}(\eta_0)$, for $1\le k \le n$, and arguing as before
we inductively prove that for $1\le k\le n$
$$
\length(\eta_k) \le \sigma^{k/2}\length(\eta_0) =
\sigma^{k/2}\dist_{f^{n}(V_n)}(f^{n}(y),f^{n}(z)) ,
$$
which implies that for $1\le k \le n$
 $$
\dist_{f^{n-k}(V_n)}(f^{n-k}(y),f^{n-k}(z)) \le
\sigma^{k/2}\dist_{f^n(V_n)}(f^{n}(y),f^{n}(z)).$$
 %If $1\le k \le n$ and $y, z\in V_n$, then we have
% \begin{align*}
%\dist_{f^{n-k}(D)}(f^{n-k}(y),f^{n-k}(z)) &\le
%\dist_{f^{n-k}(D)}(f^{n-k}(y),f^{n-k}(x)) +
%\dist_{f^{n-k}(D)}(f^{n-k}(z),f^{n-k}(x))\\
% &\le
%\sigma^{k/2}\dist_{f^n(D)}(f^{n}(y),f^{n}(x)).
% \end{align*}
This completes the proof of the first two items of the lemma.

Given $y\in V_n$ we have $\dist(f^j(x),f^j(y))\le \delta_1$ for
every $1\le j \le n$, which together with (\ref{e.delta1}) gives
\begin{eqnarray*}
% \nonumber to remove numbering (before each equation)
 % \|Df^{-k} \mid E^{cu}_{f^{n}(y)}\| &\le&
  \prod_{j=n-k+1}^{n}\|Df^{-1}
\mid E^{cu}_{f^{j}(y)}\|
  \le \sigma^{-k/2}\prod_{j=n-k+1}^{n}\|Df^{-1} \mid E^{cu}_{f^{j}(x)}\|
   \le \sigma^{k/2}.
\end{eqnarray*}
Recall that $f^j(x)\in K$ for every $j$, and $n$ is a
$\sigma$-hyperbolic time for $x$.
\end{proof}

We shall sometimes refer to the sets \( V_n \) as \emph{hyperbolic
pre-balls} and to their images \( f^{n}(V_n) \) as \emph{hyperbolic
balls}. Notice that the latter are indeed balls of radius \(
\delta_1 \).

\cco[\bf Bounded Distortion] \label{p.distortion} There exists
$C_2>1$ such that given $\Delta$ as in Lemma~\ref{l.contraction}
with $\kappa(\Delta) \le C_1$, and given any hyperbolic pre-ball
$V_n\subset \Delta$ with $n\ge n_0$, then for all $y,z\in V_n$
%
%$x\in \Delta$ and $n \ge n_0$ a $\sigma$-hyperbolic time for $x$
$$
\frac{1}{C_2} \le \frac{|\det Df^{n} \mid T_y \Delta|}
                     {|\det Df^{n} \mid T_z \Delta|}
            \le C_2.
$$
%for every $y\in \Delta$ such that $\dist_{f^n(\Delta)}(f^n(y),f^n(x))\le
%\delta_1$.
 \fco

\begin{proof}
For $0\le i <n$ and $y\in \Delta$, we denote $J_i(y)= |\det Df \mid
T_{f^i(y)} f^i(\Delta)|$. Then,
$$
\log \frac{|\det Df^{n} \mid T_y \Delta|}{|\det Df^{n} \mid T_z
\Delta|} = \sum_{i=0}^{n-1} \big(\log J_i(y) - \log J_i(z)\big).
$$
By Proposition~\ref{c.curvature}, $\log J_i$ is $(L,\zeta)$-H\"older
continuous, for some uniform constant $L>0$. Moreover, by
Lemma~\ref{l.contraction}, the sum of all $\dist_\Delta(f^j(y),
f^j(z))^\zeta$ over $0\le j \le n$ is bounded by
$\delta_1/(1-\sigma^{\zeta/2})$. Now it suffices to take
$C_2=\exp(L\delta_1/(1-\sigma^{\zeta/2}))$.
\end{proof}

%%%%%%%%%%%%%%%%%%%%%%%%%%%%%%%%%%%%%%%%%%%%%%%%%%%%%%%%%%%%%%

\section{A local unstable disk inside $\Lambda$}\label{s.local}

Now we are able to prove Theorems~\ref{t:disco1} and~\ref{t:disco2}.
These will be obtained as corollaries of the next result slightly
more general result, as we shall see next. Take $K\subset M$ a
forward invariant compact set and let $\Lambda\subset K\subset U$ be
as before.

 \ct\label{t:disco}
 Let \( f: M\to M \) be a  \( C^{1+} \)
diffeomorphism and let $K\subset M$ be a forward invariant compact
 set with a continuous splitting $T_K M=E^{cs}\oplus E^{cu}$
dominated over $\Lambda=\bigcap_{n\ge0}f^n(K)$. Assume that there is
a disk $\Delta$ tangent to the centre-unstable cone field with
$\leb_\Delta(\Delta\cap K)>0$ such that \ref{NUE} holds for every
$x\in \Delta\cap K$. Then $\Lambda$ contains some local unstable
disk.
 \ft

Let us show that  Theorem~\ref{t:disco} implies
Theorem~\ref{t:disco1}. Assume that \ref{NUE} holds for Lebesgue
almost every $x\in K$ with $\leb(K)>0$.  Choosing a $\leb$ density
point of $K$, we laminate a neighborhood of that point into disks
tangent to the centre-unstable cone field contained in $U$. Since
the relative Lebesgue measure of the intersections of these disks
with $K$ cannot be all equal to zero, we obtain some disk $\Delta$
as in the assumption of Theorem~\ref{t:disco}.

For showing that Theorem~\ref{t:disco} implies
Theorem~\ref{t:disco2}, we just have to observe that local stable
manifolds are tangent to the centre-unstable subspaces and these
vary continuously with the points in $K$, thus being tangent to the
centre-unstable cone field.

\smallskip

In the remaining of this section we shall prove
Theorem~\ref{t:disco}. Let $\Delta$ be a disk tangent to the
centre-unstable cone field intersecting $K$ in a
positive $\leb_\Delta$ subset.  %Taking a subset of $H$, if necessary,  still
%with positive $\leb_\Delta$ measure, we may assume that there is
%$c>0$ such that for every $x\in H$
%\begin{equation}
%\label{NUEc} \liminf_{n\to+\infty} \frac{1}{n}
%    \sum_{j=1}^{n} \log \|Df^{-1} \mid E^{cu}_{f^j(x)}\|\le-c.
%\end{equation}
Since \ref{NUE} remains valid under positive iteration, by
Proposition~\ref{c.curvature} we may assume that
$\kappa(\Delta)<C_1$. It is no restriction to assume that $K$
intersects the sub-disk of $\Delta$ of radius $\delta/2$, for some
$0<\delta<\delta_1$, in a positive $\leb_\Delta$ subset, and we do
so.

The following lemma is due to Pliss \cite{Pl72}, and a proof of it
in this precise form can be found in \cite[Lemma 3.1]{ABV}.

\cle %[Pliss]
\label{l.pliss} Given $A \ge c_2>c_1>0$ there exists
$\theta>0$ such that for any real numbers $a_1, \ldots, a_N$ with
$a_j \le A$ and
$$
\sum_{j=1}^{N} a_j \ge c_2 N,\quad \text{for every } 1\le j \le N,
$$
there are $l > \theta N$ and $1 < n_1 < \cdots < n_l \le N$ so that
$$
\sum_{j=n+1}^{n_i} a_j\ge c_1 (n_i-n), \quad\text{for every } 0 \le
n < n_i \text{ and  }1\le i\le l.
$$
\fle

\cco\label{l:hyperbolic2}
    There is $\sigma>0$ %depending only on $f$ and $\lambda$
     such that every \( x
    \in  \Delta\cap K \) has infinitely many $\sigma$-hyperbolic times.
     \fco
\dem Given $x\in \Delta\cap K$, by \ref{NUE} we have infinitely many
positive integers $N$ for which
$$
\sum_{j=1}^{N} \log \|Df^{-1}|E^{cu}_{f^j(x)}\| \le -\frac{c}2N.
$$
 Then it suffices to take
$c_1=c/2$, $c_2=c$, $A=\sup\big|\,\log \|Df^{-1}|E^{cu}\|\,\big|$,
and $a_j=-\log \|Df^{-1}|E^{cu}_{f^j(x)}\|$ in the previous lemma.
\cqd

Note that under assumption \ref{NUE} we are unable to prove the
existence of a positive frequency of hyperbolic times at infinity,
as in \cite[Corollary 3.2]{ABV}. This would be possible if we had
taken $\limsup$ instead of $\liminf$ in the definition of \ref{NUE}.
The existence of infinitely many hyperbolic times is enough for what
comes next.

\cle \label{l.discao} Let $O$ be an open set in $\Delta$ such that
$\leb_\Delta(O\cap K)>0$. Given any small $\rho>0$ there is a
hyperbolic time $n$, a hyperbolic pre-ball $V\subset O$ and
$W\subset V$ such that
 $\Delta_n=f^n(W)$ is a disk of radius $\delta_1/4$
tangent to the centre-unstable cone field, and $ \displaystyle
{\leb_{\Delta_n}(f^{n}(K))}\ge (1-\rho)
 {\leb_{\Delta_n}(\Delta_n)}.$
\fle

\dem Take a small number $\epsilon>0$. Let  $C$ be a compact subset
of $O\cap K$ and let $A$ be an open neighborhood of $O\cap K$ in
$\Delta$ such that
$$\leb_\Delta(A\setminus C)<\epsilon{\leb_\Delta(C)}.$$
It follows from Corollary~\ref{l:hyperbolic2} and
Lemma~\ref{l.contraction}  that we can choose for each $x\in C$ a
$\sigma$-hyperbolic time $n(x)$ and a hyperbolic pre-ball $V_x$ such
that $V_x\subset A$. Recall that $V_x$ is the neighborhood of $x$
 which is mapped diffeomorphically by
$f^{n(x)}$ onto a ball $B_{\delta_1}(f^{n(x)}(x))$ of radius
$\delta_1$ around $f^{n(x)}(x)$, tangent to the centre-unstable cone
field. Let $W_{x}\subset V_x$ be the pre-image of the ball
$B_{\delta_1/4}(f^{n(x)}(x))$ of radius $\delta_1/4$ under this
diffeomorphism. By compactness there are $x_1,...,x_m\in C$ such
that $C\subset W_{x_1}\cup...\cup W_{x_s}$. Writing
\begin{equation}\label{eqn1}
    \{n_1,...,n_s\}=\{n(x_1),...,n(x_m)\},
\quad\text{with $n_1<n_2<...<n_s$},
\end{equation}
let $I_1\subset\mathbb{N}$ be a maximal set of $\{1,...,m\}$ such
that if $i\in I_1$ then $n(x_i)=n_1$ and $W_{x_i}\cap
W_{x_j}=\emptyset$ for all $j \in I_1$ with $j\ne i$. Inductively we
define $I_k$  for $2\le k \le s$ as follows: Supposing that
$I_{k-1}$ has already been defined, let $I_k\subset\mathbb{N}$ be a
maximal set of $\{1,\dots,m\}$ such that if $i\in I_k$, then
$n(x_i)=n_k$ and $W_{x_i}\cap W_{x_j}=\emptyset$ for all $j \in I_k$
with $j\neq i$, and also $W_{x_i}\cap W_{x_j}=\emptyset$ for all
$j\in I_1\cup...\cup I_{k-1}$.

Let $I=I_1\cup\cdots\cup I_s$. By maximality, each $W_{x_j}$, for
$1\le j\le m$, intersects some $W_{x_i}$ with $i\in I$ and
$n(x_j)\ge n(x_i)$. Thus, given any $1\le j\le m$ and taking $i\in
I$ such that $W_{x_j}\cap W_{x_i}\ne\emptyset$ and $n(x_j)\ge
n(x_i)$, we get
$$
f^{n(x_i)}(W_{x_j})\cap B_{\delta_1/4}(f^{n(x_i)}(x_i))\ne\emptyset.
$$
Lemma~\ref{l.contraction}  assures that
$$
\diam(f^{n(x_i)}(W_{x_j}))\le\frac{\delta_1}{2}
\sigma^{(n(x_j)-n(x_i))/2}\le\frac{\delta_1}{2},
$$
 and so
$$f^{n(x_i)}(W_{x_j})\subset B_{\delta_1}(f^{n(x_i)}(x_i)).$$ This
implies that $W_{x_j}\subset V_{x_i}$. Hence $\{V_{x_i}\}_{i\in I}$
is a covering of $C$. It follows from Corollary~\ref{p.distortion}
that there is a uniform constant $\gamma>0$ such that
 $$
 \frac{\leb_\Delta(W_{x_i})}{\leb_\Delta(V_{x_i})}\ge\gamma,
 \quad\text{for every }i\in I.$$
Hence
 \begin{eqnarray*}
 % \nonumber to remove numbering (before each equation)
   \leb_\Delta\big(\cup_{i\in I}W_{x_i}\big)
    &=& \sum_{i\in I} \leb_\Delta(W_{x_i}) \\
    &\ge& \sum_{i\in I}\gamma \leb_\Delta(V_{x_i}) \\
   &\ge& \gamma\leb_\Delta\big(\cup_{i\in I}V_{x_i}\big)\\
   &\ge& \gamma \leb_\Delta(C).
 \end{eqnarray*}
Setting
$$
\rho=\min\left\{\frac{\leb_\Delta(W_{x_i}\setminus
C)}{\leb_\Delta(W_{x_i})}\colon i\in I\right\},
$$
we have
\begin{eqnarray*}
% \nonumber to remove numbering (before each equation)
  \varepsilon \leb_\Delta(K) &\ge& \leb_\Delta(A\setminus C) \\
   &\ge& \leb_\Delta\big(\cup_{i\in I}W_{x_i}\setminus
C\big) \\
   &\ge& \sum_{i\in I}\leb_\Delta(W_{x_i}\setminus C)\\
   &\ge& \rho \leb_\Delta\big(\cup_{i\in I}W_{x_i}\big)\\
    &\ge& \rho \gamma
\leb_\Delta(C).
\end{eqnarray*}
 This implies that $\rho<\varepsilon/\gamma$. Since $\varepsilon>0$ can
 be taken arbitrarily small, %one can find for each $\varepsilon_j=1/j$ a
%point $x_j\in K$ and a open ball of radius $\delta_1/4$
%$\Delta_j=B_{\delta_1/4}(f^{n(x_j)}(x_j))=f^{n(x_j)}(W(x_j))$ such
%that
%\begin{equation}\label{eq1}\frac{\leb_\Delta(\Delta_j\setminus
%U)}{\leb_\Delta(\Delta_j)}>1-\frac{1/\gamma^2}{j}\end{equation}
%
we may choose $W_{x_i}$ with the relative Lebesgue measure of $C$ in
$W_{x_i}$ arbitrarily close to~1. Then, by bounded distortion, the
relative Lebesgue measure of $f^{n(x_i)}(K)\supset f^{n(x_i)}(C)$ in
$f^{n(x_i)}(W_{x_i})$, which is a disk of radius $\delta_1/4$ around
$f^{n(x_i)}(x_i)$ tangent to centre-unstable cone field, is also
arbitrarily close to~$1$. Observe that since points in $K$ have
infinitely many $\sigma$-hyperbolic times, we may take the integer
$n(x_i)$ arbitrarily large, as long as $n_1$ in \eqref{eqn1} is also
taken large enough. \cqd

\cpr\label{co.seq} There are a sequence of sets $W_1\supset
W_2\supset\cdots$ and a sequence of positive integers $ n_1\le
n_2\le\cdots$ such that:
\begin{enumerate}
    \item $W_k$ is contained in some hyperbolic pre-ball with
    hyperbolic time $n_k$;
    \item  $\Delta_k=f^{n_k}(W_k)$ is a disk of
radius $\delta_1/4$, centered at some point $x_k$, tangent to the
centre-unstable cone field;
    \item $f^{n_k}(W_{k+1})$ is contained in the disk of radius
    $\delta_1/8$ centered at $x_k$;
    \item $ \displaystyle \lim_{k\to\infty}
\frac{\leb_{\Delta_{k}}(f^{n_k}(K))}{\leb_{\Delta_{k}}(\Delta_{k})}=
1
 $.
\end{enumerate}

\fpr

\dem Take a constant $0<\rho<1$ such that for any disk $D$ of radius
$\delta_1/4$ centered at some point $x$ tangent to the
centre-unstable cone field the following holds: {\em if
 $
{\leb_{D}(A)}\ge (1-\rho){\leb_{D}(D)}
 $
 for some $A\subset D$, then we must have $\leb_{D^*}(A)>0$, where $D^*\subset D$ is the
disk of radius $\delta_1/8$ centered at the same point $x$.} Note
that it is possible to make a choice of $\rho$ in these conditions
only depending on the radius of the disk and the dimension of the
disk. Surely, once we have chosen some $\rho$ satisfying the
required property, then any smaller number still has that property.

We shall use Lemma~\ref{l.discao} successively in order to define
the sequence of sets $(W_k)_k$ and hyperbolic times $(n_k)_k$
inductively. Let us start with $O=\Delta$ and $0<\rho<1$ with the
property above. By Lemma~\ref{l.discao} there are $n_1\ge1$ and
$W_{1}\subset V_1\subset O$, where $V_1$ is a hyperbolic pre-ball
with hyperbolic time $n_1$, such that $\Delta_1=f^{n_1}(W_{1})$ is a
disk of radius $\delta_1/4$ centered at some point $x_1$, tangent to
the centre-unstable cone field, with
$$
\frac{\leb_{\Delta_1}(f^{n_1}(K))}{\leb_{\Delta_1}(\Delta_1)}\ge
1-\rho
 .$$
Considering $\Delta_1^*\subset \Delta_1$ the disk of radius
$\delta_1/8$ centered at $x_1$, then by the choice of $\rho$ we have
$\leb_{\Delta_1^*}(f^{n_1}(H))>0$. Let $O_1\subset W_1$ be the part
of $W_1$ which is sent by $f^{n_1}$ diffeomorphically onto
$\Delta^*_1$. We have $\leb_\Delta(O_1\cap K)>0$.

Next we apply Lemma~\ref{l.discao} to $O=O_1$ and $\rho/2$ in the
place of $\rho$. Then we find a hyperbolic time $n_2$ and
$W_{2}\subset O_{1}$ such that $\Delta_2=f^{n_2}(W_{2})$ satisfies
  $$
\frac{\leb_{\Delta_2}(f^{n_2}(K))}{\leb_{\Delta_2}(\Delta_2)}\ge
1-\frac\rho2.$$ Observe that $W_2\subset O_1\subset W_1$. Then we
take $O_2\subset W_2$ as that part of $W_2$ which is sent by
$f^{n_2}$ diffeomorphically onto the disk $\Delta^*_2$  of radius
$\delta_1/8$ and proceed inductively. \cqd

The next proposition gives the conclusion of Theorem~\ref{t:disco}.

\cpr\label{locdisk} The sequence $(\Delta_k)_k$ has a subsequence
converging to a local unstable disk $\Delta_\infty$ of radius
$\delta_1/4$ inside~$\Lambda$. \fpr

\dem Let $(\Delta_k)_k$ be the sequence of disks given by
Proposition~\ref{co.seq} and $(x_k)_k$ be the sequence of points at
which these disks are centered. Up to taking subsequences, we may
assume that the centers of the disks converge to some point $x$.
Using Ascoli-Arzela, a subsequence of the disks converge to some
disk $\Delta_\infty$ centered at $x$. We necessarily have
$\Delta_\infty\subset \Lambda$.

Note that each $\Delta_k$ is contained in the $n_k$-iterate of
$\Delta$, which is a disk tangent to the centre-unstable cone field.
The domination property implies that the angle between $\Delta_k$
and $E^{cu}$ goes uniformly to 0 as $n\to\infty$. In particular,
$\Delta_\infty$ is tangent to $E^{cu}$ at every point in
$\Delta_\infty\subset\Lambda$. By Lemma~\ref{l.contraction}, given
any $n\ge 1$, then $f^{-n}$ is a $\sigma^{n/2}$-contraction on
$\Delta_k$ for every large $k$. Passing to the limit, we get that
$f^{-n}$ is a $\sigma^{n/2}$-contraction in the  $E^{cu}$ direction
over $\Delta_\infty$ for every $n\ge1$.
%In particular, we have shown that the subspace $E_x^{cu}$ is
%uniformly expanding for $Df$.
The fact that the $Df$-invariant splitting $T_\Lambda M=E^{cs}\oplus
E^{cu}$ is dominated implies  that any expansion $D f$ may exhibit
along the complementary direction $E^{cs}$ is weaker than the
expansion in the $E^{cu}$  direction. Then there exists a unique
unstable manifold $W^{u}_{loc}(x)$ tangent to $E^{cu}$ and which is
contracted by the negative iterates of $f$; see \cite{Pe76}. %: for every $y\in W^{uu}(x)$,
%$\dist(f^{-k}(x),f^{-k}(y))$ decreases at least as $\|Df^{-k}\mid
%E^{cu}_x\|\le \sigma^{k/2}$ when  $k$ gets large.
Since $\Delta_\infty$ is contracted by every $f^{-n}$, and  all its
negative iterates are tangent to centre-unstable cone field, then
$\Delta_\infty$ is contained in $W^{u}_{loc}(x)$.
%Note that every point in $f^{k_n}(H)$ still has non-uniform
%expansion in the $E^{cu}$ direction, since the set of those points
%is forward invariant, by definition. Using the absolute continuity
%of the strong stable foliation and the fact that if $z$ has
%non-uniform expansion in the $E^{cu}$ direction and $w\in W^s(z)$
%then $w$ also has non-uniform expansion in the $E^{cu}$ direction,
%we get that $ m_D$ almost every point in $D$ has non-uniform
%expansion in the $E^{cu}$.
 \cqd

%%%%%%%%%%%%%%%%%%%%%%%%%%%%%%%%%%%%%%%%%%%%%%%%%%%%%%%%%%%%%%%%%%%%%%%

\section{Existence of hyperbolic periodic  points}\label{s.hypeper}

Here we prove Theorem~\ref{t:limitph}. By Proposition~\ref{co.seq}
there exist a sequence of sets $W_1\supset W_2\supset\cdots$
contained in $\Delta$ and a sequence of positive integers $ n_1\le
n_2\le\cdots$ such that:
\begin{enumerate}
    \item $W_k$ is contained in some hyperbolic pre-ball with
    hyperbolic time~$n_k$;
    \item  $\Delta_k=f^{n_k}(W_k)$ is a disk of
radius $\delta_1/4$, centered at some point $x_k$, tangent to the
centre-unstable cone field;
    \item $f^{n_k}(W_{k+1})$ is contained in the disk $\Delta_k^*$ of radius
    $\delta_1/8$ centered at $x_k$.
\end{enumerate}
 Taking a subsequence, if necessary,
we have by Proposition~\ref{locdisk} that the sequence of disks
$(\Delta_k)_k$ accumulates on a local unstable disk $\Delta_\infty$
of radius $\delta_1/4$ which is contained in~$\Lambda$.
 Our aim now is to prove that
$\Lambda$ contains the unstable manifold of some periodic point.

Similarly to \eqref{e.delta1}, we choose $\delta>0$ small so that
$W^s_\delta(z)$ is defined for every $z\in\Lambda$, the
$2\delta$-neighborhood of $\Lambda$ is contained in $U$, and
\begin{equation}
\label{e.delta} \|Df^{-1}(f(y)) v \| \le \sigma^{-1/4} %\frac{1}{\sqrt[4]\sigma}
\|Df^{-1}| E^{cu}_{f(x)}\|\,\|v\|,
\end{equation}
whenever $x\in U$, $\dist(x,y)\le 2\delta$, and $v\in C^{cu}_a(y)$.

%\cle\label{l.fibra} There are a point $z\in \Delta_\infty$ and an
%integer $n\ge 1$ such that $f^n(W^s_{\delta}(z))\subset
%W^s_{\delta}(z)$. Moreover, we may have $\dist(f^n(W^s_{\delta}(z),
%z)$ arbitrarily small.\fle

\cpr\label{pr.p} Given $\Lambda_1\subset \Lambda$ with
$\leb(\Lambda_1)>0$, there exist a hyperbolic periodic point
$p\in\Lambda$ and $\delta_2>0$ (not depending on $p$) such that:
\begin{enumerate}
    \item $\overline{W^u(p)}\subset \Lambda$;
    \item the size of ${W^u_{\loc}(p)}$ is at least $\delta_2$;
    \item $\leb_{W^u_{\loc}(p)}$ almost every point in $W^u_{\loc}(p)$
belongs to $H$;
 \item there is $x\in \Lambda_1$ with $\omega(x)\subset
 \overline{W^u(p)}$.
\end{enumerate}

  \fpr

\dem Let $x$ denote the center of the accumulation disk
$\Delta_\infty$. Let us consider the cylinder
 $$\cc_\delta=\bigcup_{y\in \Delta_\infty}W_\delta^s(x),$$
 and  the
projection along local stable manifolds
$$
\pi\colon\cc_\delta\longrightarrow\Delta_\infty.
$$
Slightly diminishing  the radius of the disk $\Delta_\infty$, if
necessary, we may assume that there is a positive integer $k_0$ such
that for every $k\ge k_0$
\begin{equation}\label{prje}
    \pi(\Delta_k\cap\cc_\delta)=\Delta_\infty\qand \Delta_k^*\subset\cc_\delta.
\end{equation}
For each $k\ge k_0$ let
$$\pi_k\colon\Delta_\infty\longrightarrow\Delta_k$$ be the
projection along the local stable manifolds. Notice that these
projections are continuous and
$\pi\circ\pi_k=\text{id}_{\Delta_\infty}$. Take a positive integer
$k_1> k_0$ sufficiently large so that
 \begin{equation}\label{ka1}
    \pi(\Delta_{k_1}\cap
 \cc_{\delta/2})=\Delta_\infty\qand \lambda^{n_{k_1}-n_{k_0}}\le\frac14.
 \end{equation}
We have
 $$\Delta_{k_1}=f^{n_{k_1}}(W_{k_1})\subset
 f^{n_{k_1}-n_{k_0}}\left(f^{n_{k_0}}(W_{k_0+1})\right)\subset
 f^{n_{k_1}-n_{k_0}}(\Delta_{k_0}^*),
 $$
which together with~\eqref{prje} and \eqref{ka1} implies that there
is some disk $\Delta_0\subset\Delta_\infty$ such that
 $$
 \pi\circ
 f^{n_{k_1}-n_{k_0}}\circ\pi_{k_0}(\Delta_0)=\Delta_\infty.$$
Thus there must be some $z\in\Delta_0\subset\Delta_\infty$ which is
a fixed point for the continuous map $\pi\circ
 f^{n_{k_1}-n_{k_0}}\circ\pi_{k_0}$. This means that there are
 $z_{k_0},z_{k_1}\in W^s_\delta(z)$ with $z_{k_0}\in\Delta_{k_0}$
 and $z_{k_1}\in\Delta_{k_1}$ such that
 $f^{n_{k_1}-n_{k_0}}(z_{k_0})=z_{k_1}$. Letting
 $\gamma=W^s_\delta(z)$, we have   $\dist_\gamma(w,z_{k_1})\le2\delta$
 for every $w\in\gamma$.
This implies that
 $$\dist_\gamma(f^{n_{k_1}-n_{k_0}}(w),z_{k_1})=
  \dist_\gamma(f^{n_{k_1}-n_{k_0}}(w),f^{n_{k_1}-n_{k_0}}(z_{k_0}))\le2\delta\lambda^{n_{k_1}-n_{k_0}},$$
which together with \eqref{ka1} gives
$$
 \dist_\gamma(f^{n_{k_1}-n_{k_0}}(w),z)\le
 \dist_\gamma(f^{n_{k_1}-n_{k_0}}(w),z_{k_1})+\dist_\gamma(z_{k_1},z)\le
 \delta.
 $$
 We conclude that $f^{n_{k_1}-n_{k_0}}(W^s_\delta(z))\subset
 W^s_\delta(z)$.
% Moreover, since we may take $\Delta_{k_1}$ arbitrarily close to
% $\Delta_\infty$, the second conclusion
% of the lemma also follows  if we take $k_1$ sufficiently large.\cqd
% \cre
%The proof of Lemma~\ref{l.fibra} gives that
% \fre
Since $W^s_{\delta}(z)$ is a topological disk, this implies that
$W^s_{\delta}(z)$ must necessarily contain some periodic point $p$
of period $m=n_{k_1}-n_{k_0}$. As $z\in\Delta_\infty$ and $p\in
W^s_{\delta}(z)$ it follows that $p\in\Lambda$, by closeness of
$\Lambda$.

Let us now prove that $p$ is a hyperbolic point. As
 $p\in W^s_{\delta}(z)$, it is enough to show that
 $\|Df^{-m} \mid E^{cu}_{f^m(p)}\| <1.$
 Let $q=W^s_\delta(z)\cap f^{n_{k_0}}(W_{k_1})$. Observe that since
$p\in\Lambda\cap W^s_\delta(z)$, then $q$ belongs to the
$2\delta$-neighborhood of $\Lambda$, which is contained in $U$.
Since $W_{k_1}$ is contained in some hyperbolic pre-ball with
hyperbolic time~$n_{k_1}$, it follows from Lemma~\ref{l.contraction}
that for every $1\le j \le n_{k_1}$ and $y\in W_{k_1}$,
$$
\|Df^{-j} \mid E^{cu}_{f^{n_{k_1}}(y)}\| \le  \sigma^{j/2}.
$$
In particular, taking $j=m=n_{k_1}-n_{k_0}$ and $y=f^{-n_{k_0}}(q)$,
we have
 $$
\|Df^{-m} \mid E^{cu}_{f^{m}(q)}\| \le  \sigma^{m/2}.
 $$
The choice of $\delta$ in \eqref{e.delta} together with the fact
that $p,q\in W^s_\delta(z)$ imply that
\begin{eqnarray}
% \nonumber to remove numbering (before each equation)
  \|Df^{-m} \mid E^{cu}_{f^{m}(p)}\| &\le& \prod_{j=1}^{m}\|Df^{-1}
\mid E^{cu}_{f^{j}(p)}\| \label{pontyp}\\
  &\le& \sigma^{-m/4}\prod_{j=1}^{m}\|Df^{-1} \mid E^{cu}_{f^{j}(q)}\| \nonumber\\
   &\le& \sigma^{m/4}.\label{contsigma}
\end{eqnarray}
Thus we have proved the hyperbolicity of $p$.

Now since $p$ is a hyperbolic periodic point, there is
$W^u_{\loc}(p)$ a local unstable manifold through $p$ tangent to the
center unstable bundle. As $\Delta_\infty$ cuts transversely the
local stable manifold through $p$, then using the $\lambda$-lemma we
deduce that the positive iterates of $\Delta_\infty$ accumulate on
the unstable manifold through $p$. Since these iterates are all
contained in $\Lambda$ and $\Lambda$ is a closed set, we must have
$W^u(p)\subset\Lambda$, which then implies that
$\overline{W^u(p)}\subset\Lambda$. Thus we have proved the first
part of the result.

By \eqref{pontyp} and \eqref{contsigma} we deduce that every multiple of $m$ is a
$\sigma^{1/4}$-hyperbolic time for $p$. Then we choose $\delta_2>0$ such that
an inequality as in \eqref{e.delta1} holds with $\delta_2$ in the place of
$\delta_1$ and $\sigma^{1/8}$ in the place of $\sigma^{1/2}$.
Using Lemma~\ref{l.contraction} with $W^u_\loc(p)$ in
the place of $\Delta$ and taking a sufficiently large $\sigma^{1/4}$-hyperbolic time for $p$
we deduce that there is a hyperbolic pre-ball inside $W^u_\loc(p)$. This imples
that its image by the hypebolic time, which is a disk of radius $\delta_2$
around $p$, is contained in the local unstable manifold of $p$. This
gives the second part of the result.

 Observe that as long as we take the local
unstable manifold through $p$ small enough, then every point in $W^u_{\loc}(p)$
belongs to the local stable manifold of some point in
$\Delta_\infty$.
By construction, $\Delta_\infty$ is accumulated by
the disks $\Delta_k=f^{n_k}(W_k)$ which, by
Proposition~\ref{co.seq}, satisfy
 \begin{equation}\label{fnkh}
  \lim_{k\to\infty}
\frac{\leb_{\Delta_{k}}(f^{n_k}(H))}{\leb_{\Delta_{k}}(\Delta_{k})}=
1.
 \end{equation}
Since $H$ is positively invariant, we have
 $$
 \lim_{k\to\infty}
\frac{\leb_{\Delta_{k}}(H)}{\leb_{\Delta_{k}}(\Delta_{k})}= 1.
$$ Let now
$\varphi\colon\Lambda\to\RR$ be the continuous function given by
 $$\varphi(x)=\log \|Df^{-1} \mid E^{cu}_{x}\|.$$
Since Birkhoff's time averages are constant for points in a same
local stable manifold and the local stable foliation is absolutely
continuous, we deduce that
 $$
\frac{\leb_{\Delta_{\infty}}(H)}{\leb_{\Delta_{\infty}}(\Delta_{\infty})}=
1.
$$
The same conclusion holds  for the local unstable manifold of $p$ in
the place of $\Delta_{\infty}$ by the same reason.

Let us now prove the last item. Since $H$ has full Lebesgue measure
in $\Lambda$ and $\Lambda_1\subset\Lambda$ has positive Lebesgue
measure, we may start our construction with the set $H_1=H\cap
\Lambda_1$ in the place of $H$ intersecting the disk $\Delta$ in a
positive $\leb_\Delta$ measure set of points. Although we have not
invariance of $H_1$, by \eqref{fnkh} we still have the property that
the iterates of $H_1\subset\Lambda_1$ accumulate on the whole
$\Delta_\infty$. Since the stable manifolds through points in
$W^u_{\loc}(p)$ intersect $\Delta_\infty$, there must be points in
$\Lambda_1$ accumulating on $W^u_{\loc}(p)$.
 \cqd

Let $p_1$ be a hyperbolic periodic point as in
Proposition~\ref{pr.p}. Let $B_1$ be the basin of
$\overline{W^u(p_1)}$, i.e. the set of points $x$ whose
$\omega$-limit is contained in $\overline{W^u(p_1)}$. If
$\leb(\Lambda\setminus B_1)=0$, then we have proved the theorem.
Otherwise, let $\Lambda_1=\Lambda\setminus B_1$. Using again
Proposition~\ref{pr.p} we obtain a point $p_2\in\Lambda$ such that
the basin $B_2$ of $\overline{W^u(p_2)}$ attracts some point of
$\Lambda_1$. By definition of $\Lambda_1$ we must have
$\overline{W^u(p_1)}\neq \overline{W^u(p_2)}$.

We proceed inductively, thus obtaining periodic points $p_1,\dots,
p_n\in\Lambda$ with $\overline{W^u(p_i)}\neq \overline{W^u(p_j)}$
for every $i\neq j$. This process must stop after a finite number of
steps. Actually, if there were an infinite sequence of points as
above, by compactness, choosing $p_{i_1},p_{i_2}$ sufficiently
close, using the inclination lemma we would get
$\overline{W^u(p_{i_1})}= \overline{W^u(p_{i_2})}$.

So far we have proved the first two items of
Theorem~\ref{t:limitph}. Assume now that $E^{cu}$ has dimension one.
We want to show that each $\overline{W^u(p_i)}$ attracts an open set
containing $ \overline{W^u(p_i)}$. Given $1\le i\le k$, by
Proposition~\ref{pr.p} we can find at least one point on each
connected component of $W^u(p_i)\setminus \{p_i\}$ belonging to $H$.
Since these these points have infinitely many hyperbolic times, then
each connected component of $W^u(p_i)\setminus \{p_i\}$ must
necessarily have infinite arc length; recall
Lemma~\ref{l.contraction}. This implies that each point $x\in
W^u(p_i)$ has an unstable arc $\gamma^u(x)\subset W^u(p_i)$ of a
fixed length passing through it. Let
 $$B(x)=\bigcup_{y\in\gamma^u(x)}W^s_\delta(y).$$
By domination, the angles of $\gamma^u(x)$ and the local stable
manifolds $W^s_\delta(y)$ with $y\in \gamma^u(x)$ are uniformly
bounded away from zero. Thus, $B(x)$ must contain some ball of
uniform radius (not depending on $x$), and so the set $\bigcup_{x\in
W^u(p_i)} B(x)$ is a neighborhood of $ \overline{W^u(p_i)}$. Since,
for each $x\in W^u(p_i)$, the points in $B(x)$ have their
$\omega$-limit set contained in $ \overline{W^u(p_i)}$, we are done.

%%%%%%%%%%%%%%%%%%%%%%%%%%%%%%%%%%%%%%%%%%%%%%%%%%%%%%%%%%%%%%%%%%%%%%%%

\section{Hyperbolic sets with positive volume}\label{s.hyperbolic}

In this section we prove Theorem~\ref{t:Anosov} and
Theorem~\ref{t:limit}. %, and give an example which illustrates the
%necessity of the transitivity assumption in Theorem~\ref{t:Anosov}.
 Since in the present
situation $Df\vert E^u_x$ is uniformly expanding, then we
have~\ref{NUE}  for every $x\in\Lambda$.

\subsection{Transitive case} Assume first that $\Lambda$ has
positive volume. It follows from Corollary~\ref{c:disco1} that
$\Lambda$ must contain some local unstable disk. The first item of
Theorem~\ref{t:Anosov} is a consequence of the following folklore
lemma whose proof we give here for the sake of completeness.

\cle\label{l.unstable} If $\Lambda$ is a transitive hyperbolic set
containing the local unstable manifold of some point, then $\Lambda$
contains the local unstable manifolds of all its points. \fle
 \dem
 Take
$\delta>0$ small such that $W^s_{\delta}(x)$ and $W^u_{\delta}(y)$
intersects at most in one point, for every $x,y\in\Lambda$, and
assume that $W^u_{\delta}(x_0)\in\Lambda$ for some $x_0\in\Lambda$.
Let $z\in \Lambda$ be a point with dense orbit in $\Lambda$. It is
no restriction to assume that $W^s_{\delta}(z)$ intersects
$W^u_{\delta}(x_0)$, and we do so. Let $x_1=W^s_{\delta}(z)\cap
W^u_{\delta}(x_0)$. We also have $W^u_{\delta}(x_1)\subset\Lambda$.
Given any point $y\in\Lambda$, we take a sequence of integers
$0=n_1<n_2<\cdots$ such that $f^{n_k}(z)\to y$, when $k\to\infty$.
Since $x_1\in W^s(z)$ we also have $x_k:=f^{n_k}(x_1)\to y$, when
$k\to\infty$. The local unstable manifolds through the points
$x_1,x_2,...$ are necessarily contained in $\Lambda$ and accumulate
on a disk $D(y)$ contained in $\Lambda$ and containing $y$. Since
the local unstable disks are tangent to the unstable spaces, the
continuity of these spaces implies that $T_wD(y)=E_w^u$ for every
$w\in D(y)$. By uniqueness of the unstable foliation, we must have
$D(y)$ coinciding with the local unstable manifold through $y$. \cqd

Using the previous lemma applied to $f^{-1}$, we have that $\Lambda$
must also contain the stable manifolds through its points. Then we
easily deduce that every point in $\Lambda$ belongs in the interior
of $\Lambda$, thus showing that $\Lambda$ is an open set. Since
$\Lambda $ is assumed to be closed, we conclude that $\Lambda=M$,
thus having proved the first part of Theorem~\ref{t:Anosov}.

\cle\label{l.point} Let $\Lambda$ be a hyperbolic set attracting a
set with positive volume. Then there is a point in $\Lambda$ whose
local unstable manifold is contained in~$\Lambda$. \fle

\dem We fix continuous extensions (not necessarily continuous) of
the two bundles $E^{cs}$ and $E^{cu}$ to some neighborhood $U$ of
$\Lambda$. Let $A$ be the set of points which are attracted to
$\Lambda$ under positive iteration. Since $A$ has positive volume,
there must be some compact set  $C \subset A$ with positive volume,
and some $N\in\NN$ such that $f^n(C)\subset U$ for every $n\ge N$.
 Letting
  $$K=\bigcup_{n\ge N} f^n(C) \cup \Lambda$$ we have that $K$ is compact forward invariant set with
positive volume for which $$\Lambda=\bigcap_{n\ge1}f^n(K).$$ The
conclusion of the lemma then follows from Theorem~\ref{t:limitph}.
\cqd

The second part of Theorem~\ref{t:Anosov} can be now easily deduced
from Lemma~\ref{l.unstable} and Lemma~\ref{l.point}. Actually, it
follows from the lemmas that
 $$\bigcup_{x\in\Lambda}W^s_\delta(x)$$
is a neighborhood of $\Lambda$ whose points are attracted to
 $\Lambda$ under positive iteration.

\subsection{Nontransitive case} Here we consider the case hyperbolic sets with
positive volume not necessarily transitive and prove
Theorem~\ref{t:limit}.

Let $\Sigma=\overline{W^u(p)}\subset\Lambda$, where $p$ is a
hyperbolic periodic point given by Proposition~\ref{pr.p}. We claim
that $\Sigma$ contains the local unstable manifolds of all its
points. Indeed, if $x\in\Sigma$, then there is a sequence $(x_n)_n$
of points in $W^u(p)$ converging to~$x$. The continuous variation of
the local unstable manifolds gives that the local unstable manifolds
of the points $x_n$, which are contained in $\Sigma$, accumulate on
the local unstable manifold of $x$. By closeness, the local unstable
manifold of $x$ must  be contained in $\Sigma$. Thus, defining
 $$A=\bigcup_{x\in\Sigma}W^s_\delta(x)$$
we have that $A$ is a neighborhood of $\Sigma$ whose points have
their $\omega$-limit set contained in $\Sigma$. Since $\Sigma$ is a
hyperbolic set with a local product structure attracting an open
neighborhood of itself, then by \cite[Theorem 18.3.1]{KH} there are
hyperbolic invariant sets
$\Omega_1,\dots,\Omega_s\subset\Sigma\subset\Lambda$ verifying
(3)-(6) of Theorem~\ref{t:limit}. Moreover, their union is
 the set of non-wandering points of $f$ in~$\Sigma$,
 $$NW(f\vert\Sigma)=\Omega_1\cup\cdots\cup\Omega_s.$$
Since $L(f\vert\Sigma)\subset NW(f\vert\Sigma)$, this implies that
$\omega(x)\subset \Omega_1\cup\cdots\cup\Omega_s$ for every $x\in
A$. Recall that every point in $A$ belongs to the stable manifold of
some point in $\Sigma$. Now since  $\Omega_1,\dots,\Omega_s$ are
disjoint compact invariant sets, given $x\in A$, we must even have
$\omega(x)\subset \Omega_i$ for some $1\le i\le s$. Reordering these
sets if necessary, let $\Omega_1,\dots,\Omega_q$, for some $q\le s$,
be those which attract a set with positive Lebesgue measure. By
Theorem~\ref{t:Anosov} and transitivity, each
$\Omega_1,\dots,\Omega_q$ attracts a neighborhood of itself.

\end{document}